# Mathematical Modeling of Sediments in the Filter and Improvement of the Filter Construction

Yuri Vitalievich Troshchiev[a]

[a]Dept. of Calculate Mathematics and Cybernetics, Moscow State University, Moscow, Russia; Ph.D., Senior Researcher; E-mail: yuvt@yandex.ru.

**Abstract.**

The filters work in many areas of technology. There constructions are different and substances under filtration are different. It is necessary in some cases to take into account forming of sediments on the walls of the filter since they can change properties of the filter or blind the filtering apertures at all. I construct a mathematical model of sedimentation growth on the walls of the porous filter in this article. Analytical investigation is present in the article and numeric results too. There are formulas for dependencies of concentration near the walls on inner concentration in liquid in the article. Flow speed, calculated time of work, purification efficiency and other parameters proved to be important factors. Differing of radiuses of apertures from membrane to membrane can make contamination equal along the filter. Numerical results show importance of preliminary calculation of the filter for the purpose it will serve. Forming of calcic sediment is an investigated example of chemical reaction.

**Introduction.**

Modern approach to development of technical devices implies high-technological constructing. Concerning the filters, it is usually possible to improve such a construction as simply porous material that retains small particles. Many investigators have done rather much in this direction, for example, [1] [2] [3] [4] [5] [6] [7] [8] [9] [10] [11] [12] [13] [14]. The majority of these articles are reviews and represent diverse enough examples both about constructions of the filters and about mixtures that pass through the filters. For example, filtration of polysaccharides [10], bacteria [12] or microplastics [14].

To simplify terminology let the matter that the filter retains be contamination in spite of the fact that sometimes contamination is the matter which the filter does not retain or the filter separates two useful matters.

The author proposed stochastic mathematical model of the porous filter which supposes ordered arrangement of the cells and different size of the filtering apertures [15] [16] [17]. It is possible to produce such a filter, for example, by 3D-printer or membranes and laser formation of apertures. It is possible to investigate the properties of such filters analytically and by computer simulation too. There were some interesting results. For example, dirt accumulates in adjacent cells in crosscut layers. Another one: the filter assembled of the membranes is more rational, i.e. it is useful to increase apertures in the sides of the cells to make free flow of liquid in the crosscut directions. It is also interesting that contamination from layer to layer gather in something like chess order, i.e. the contaminated area follows the free area in the previous membrane and vice versa.

Concerning this article, I construct a mathematical model of sediment growth on the inner surface of the filter cells and investigate interaction of filtration and sedimentation. I.e. one more process is included in mathematical model of the filter additionally to filtration. Sediments (for instance, calcic sediment) are generally considerable problem and they can spoil functionality of high-technological filter. The article also represents an example of chemical reaction that can cause growth of sediment.

Other examples of layers on the membrane surfaces are in [2] and [10].

**Constructing Mathematical Model of Forming of Sediments.**

Let us consider filtering aperture of radius $R$. It is worth to note that sedimentation is not the inverse process for dissolution. Indeed, if solution of the sedimentation substance is not saturated then dissolving takes place. Inverse process can be only in such a case as supersaturated solution. Thus, chemical reactions proceed while sedimentation. Let there is already an initial layer of the sedimentation substance on the walls because its sedimentation on the clean walls will proceed by another chemical scheme.

Under assumption of steady flow the dependence of velocity of fluid on radius is

$$v(r) = v_0\left(1 - \frac{r^2}{R^2}\right), \tag{1}$$

where $r$ – distance from the center of the aperture, $v_0$ – velocity in the center of the aperture. Thus, there is almost immoveable fluid layer near the wall where diffusion plays the main role [18] [19] [3] [11]. Let velocity of sedimentation is determined by concentration of some dissolved substance [11], i.e. the chemical reaction at the surface is fast (variant of such a reaction is in what follows). Let entrance concentration of this substance (may be in dissociated form) in the aperture is $c_0$ and concentration near the wall is $c_1$. Velocity of sedimentation reaction can be written as $Kc_1^n$, where $n$ – number of molecules or ions of dissolved substance taking part in the single reaction of appearing of the sediment molecule (molecules). Let the boundary of the slow layer is such that diffusion flow equals convection flow there. Then we obtain a system of equations:

$$\begin{aligned} v &= v_0\left(1 - \frac{r^2}{R^2}\right), \\ Kc_1^n &= c_0 v, \\ D(c_0 - c_1) &= c_0 v(R - r), \end{aligned} \tag{2}$$

where $D$ – diffusion coefficient, $r \in (0, R)$ – distance from the center to the boundary of the slow layer, $c_1 \in (0, c_0)$.

Introducing notation $y = 1 - r/R \in (0,1)$, $a = KR/D > 0$ и $b = K/v_0 > 0$, we obtain at $n = 1$ an equation

$$y(2-y)(ay+1) - b = 0, y \in (0,1). \tag{3}$$

For $r$ be positive it is necessary and sufficient the following inequality

$$v_0 > v_{stat} = \frac{K}{1 + KR/D}. \tag{4}$$

In the pointed out interval $y$ is increasing function of $b$ and, correspondingly, decreasing function of $v_0$. Additionally it is possible to show that if fluid flows in a pipe of radius $R$ with velocity $v_0(R) = v_{stat}(R)$, and then the pipe becomes narrower to radius $r < R$, then in the more narrow part of the pipe the following inequality is true

$$v_0(r) = v_0(R)R^2/r^2 > v_{stat}(r) = K/(1 + Kr/D). \tag{5}$$

Concerning any natural $n$ we obtain instead of (3) the following equation

$$Kc_1^n = c_0 v_0 \frac{D(c_0 - c_1)}{Kc_1^n R}\left(2 - \frac{D(c_0 - c_1)}{Kc_1^n R}\right). \tag{6}$$

It is convenient to use this equation to find $c_1$ at any $n$ and at $n = 1$ as well.

Let us denote $f_1 = D(c_0 - c_1)/(Kc_1^n R)$, $f_2 = f_1(2 - f_1)$, $f = c_0 v_0 f_2$. If $c_1$ increases from zero to $c_0$ then the left part of the equation increases from zero to $Kc_0^n$. The function $f_1$ decreases from infinity to zero when $c_1$ increases from zero to $c_0$. Additionally, this function is $y$ and must be in the interval $(0,1)$ according to physical sense. The function $f_2 = y(2 - y)$ increases from zero to unit at this interval. Thus,

$$Kc_1^n < c_0 v_0 \tag{7}$$

is a necessary and sufficient condition for existence of the physically sensed solution of the equation (6). At $n = 2$ we obtain

$$v_0 > v_{stat} = Kc_0^{-1}\left(\frac{-D + \sqrt{D^2 + 4KRDc_0}}{2KR}\right)^2. \tag{8}$$

If inequality $v_0 > v_{stat}$ ((4), (7), (8)) is violated then we obtain the following equation instead of the system of equations (2)

$$D(c_0 - c_1) = Kc_1^n R, \tag{9}$$

where $c_0$ – concentration in the center of the aperture.

**Velocity of sediment growth.**

Let the wall between the cells is thin, i.e. the length of the aperture is small. Then it is possible to neglect decrease of concentration along the aperture.

Velocity of sediment growth is determined by the value $Kc_1^n$:

$$\frac{ds}{dt} = \frac{Kc_1^n \mu_2 n_2}{n \rho_2 N_A}, \tag{10}$$

where $s$ – width of the sediment layer, $\mu_2$ – molar mass of the sediment substance, $n_2$ – number of molecules of the sediment appearing in the single chemical reaction, $\rho_2$ – density of the sediment substance, $N_A$ – Avogadro constant. In the case of the system (2) at $n = 1$

$$c_1 = \frac{c_0 v_0 y(2-y)}{K} = \frac{c_0}{ay+1}. \tag{11}$$

In the case of the equation (9) at $n = 1$

$$c_1 = \frac{Dc_0}{(KR + D)}, \tag{12}$$

and at $n = 2$

$$c_1 = \frac{-D + \sqrt{D^2 + 4KRDc_0}}{2KR}. \tag{13}$$

Concentration $c_1$ is an increasing function of $v_0$ in the formula (11) and in the formula (12) it does not depend on $v_0$.

**A variant of chemical reaction.**

The most common sediment is probably calcic one. For example, some quantity of calcium sulfate $CaSO_4$ almost always exists in water. Sodium carbonate $Na_2CO_3$ also can be present. These substances react producing almost insoluble calcium carbonate $CaCO_3$ which is the main component of calcic sediment. In this case, calcium carbonate appears according to the following chemical equation

$$CaSO_4 + Na_2CO_3 =\downarrow CaCO_3 + Na_2SO_4. \tag{14}$$

Sodium carbonate is soapy a little and adheres to the walls. Thus, the walls of the aperture will be covered by the layer of sodium carbonate, and concentration of calcium ions will control the reaction velocity. Some quantity of calcium carbonate also appears in water far away from the wall, certainly.

Let, for example, calcium carbonate concentration is $0.001\ kg/m^3$, water flows slowly and calcic sediment grows $0.1\ mm$ per month. Then it follows from (10) and (12)

$$c_0 = \frac{c_{0,m} N_A}{\mu_0},$$

$$c_1 = \frac{Dc_0}{KR + D}, \tag{15}$$

$$K = \frac{s' n \rho_2 \mu_0 D}{Dc_{0,m}\mu_2 - s' n \rho_2 \mu_0 R}.$$

Here $c_{0,m}$ – mass concentration of the dissolved substance ($CaSO_4$), $\mu_0$ – molar mass of the dissolved substance, $s' = ds/dt$ – velocity of sediment growth. It is necessary and sufficient for $K$ be greater than zero to satisfy the inequality

$$Dc_{0,m}\mu_2 - s'n\rho_2\mu_0 R > 0. \qquad (16)$$

If $D = 10^{-9} \ m^2/s$, a $R = 10^{-6} \ m$ then $K = 1.66 \cdot 10^{-4} \ m \ s^{-1}$, $c_0 = 4.4 \cdot 10^{21} \ m^{-3}$, $c_1 = 3.8 \cdot 10^{21} \ m^{-3}$. The minimal value of velocity from formula (4) is $v_0 > 1.4 \cdot 10^{-4} \ m/s$.

The more usual concentration is $0.2 \ kg/m^3$ and more for natural water. Additionally, it is possible to increase flow velocity. Thus, it is interesting in what cases it is possible to neglect concentration decrease along the filter.

**Decrease of sediment substance concentration along the filter.**

It is impossible to neglect decrease of concentration at zero flow speed because after a long enough time all the substance will precipitate. And if the speed tends to infinity then concentration near the wall tends to $c_0$ (11), the flow of particles to the wall because of chemical reaction tends, so, to $Kc_0 < \infty$. As time of liquid passing through the filter tends to zero, so, it is possible to neglect decrease of concentration if speed tends to infinity.

Let us determine firstly when it is possible to neglect the decrease if $v_0 \leq v_{stat}$. We neglect change of concentration and seek when it leads to contradiction. Change of concentration at pipe segment of the length $L$ is

$$\Delta c_0 \approx 4\pi R \frac{Kc_1^n L}{v_0}, \qquad (17)$$

where $c_1$ is defined by formula (12) at $n = 1$ and (13) at $n = 2$. So, it is really possible to neglect change of concentration if $\Delta c_0 \ll \langle c \rangle = c_0/3 + 2c_1/3$. The radius $R$ here is not the radius of aperture but typical size of cells, i.e. of cavities connected by filtering apertures, because the most mass of sediment will be there. Average concentration $\langle c \rangle$ is of magnitude of $c_0$, so, it is possible to compare $\Delta c_0$ and $c_0$: $\Delta c_0 \ll c_0$. It is always possible to choose so small $v_0$ that this inequality will be false.

Let now $v_0 > v_{stat}$. Formula (17) determines change of concentration but $c_1$ is from formula (11) or from equation (6) now. The following more complicated formula determines average concentration

$$\langle c \rangle = \frac{c_0}{R^2}\left(r^2 + R(R+r) - \frac{2(R^3 - r^3)}{3(R-r)}\right) + \\ + \frac{c_1}{R^2}\left(-r(R+r) + \frac{2(R^3 - r^3)}{3(R-r)}\right), \quad (18)$$

but it majorizes $\langle c \rangle$ for slow speed. So, it is also possible to check the inequality: $\Delta c_0 \ll c_0$.

**Construction of the filter.**

Let the filter consists of the ordered cells connected by filtering apertures (Fig. 1). Two sets of the cells on the opposite sides of the filter form inlet and outlet. External facets of these cells are open. External facets of all the other cells are blindly closed. The apertures are in every inner facet of each cell. There radiuses can differ. Flowing through the aperture the particle can block it, i.e. the particle stops and flow through the aperture stops. Sizes of the filter along the axes $L_x, L_y, L_z$, number of the cells $n_x, n_y, n_z$ and sizes of the cells, so, $h_x = L_x/n_x$, $h_y = L_y/n_y$, $h_z = L_z/n_z$.

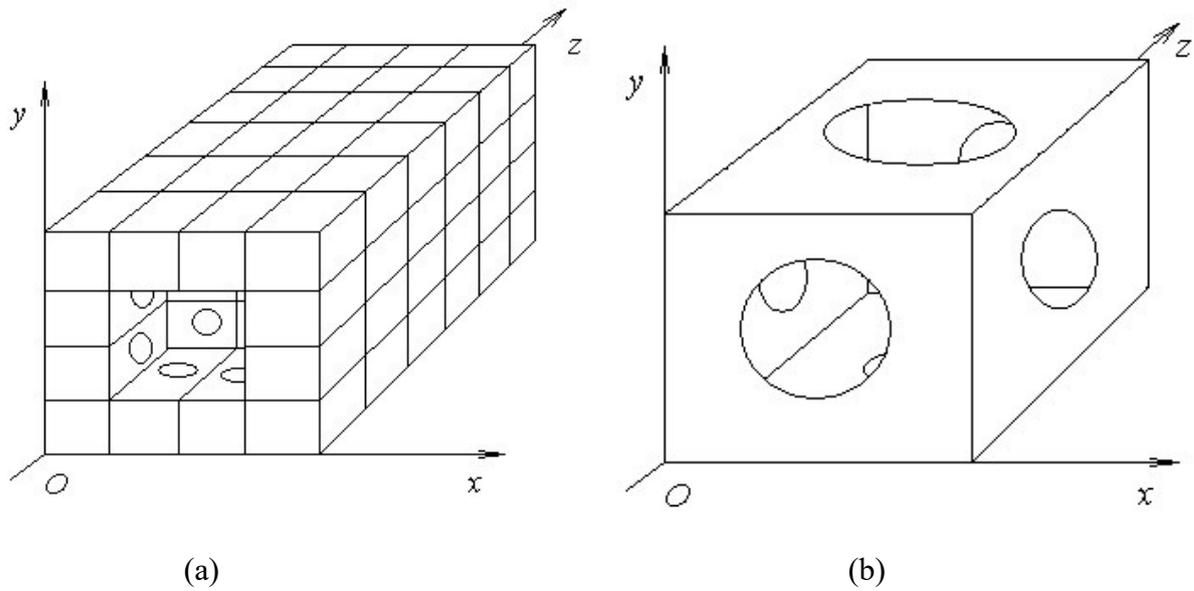

(a)             (b)

*Fig. 1 a) Filter. b) Cell.*

The size of apertures can differ now during the filter work because of sediment growth. Thus, we need formula to calculate the probability of blocking depending on the radius $r$ of the aperture. Let the particles are thin rods of the length $l$ and if such a particle comes to the aperture so that its projection is longer that the diameter of the aperture $d = 2r$ then it blocks the aperture. If the length $l$ is less than $2r$ then the probability of infiltration is 1. Let now $l \geq 2r$. Then it is possible to consider a sphere of radius $l$ and two identical disks on its surface on the opposite sides (similar to Northern and Southern polar circles on the Earth). Projections of the disks in direction of the sphere diameter connecting their centers are of radius $r$. The probability of infiltration is the sum of squares of the disks divided by square of all the sphere:

$$q = \frac{2\pi l\left(l/2 - \sqrt{l^2/4 - r^2}\right)}{\pi l^2} = 1 - \sqrt{1 - \left(\frac{2r}{l}\right)^2}. \tag{19}$$

Despite this formula provides some inaccuracy, the principle itself is applicable to obtain formulas for different objects, including filtration of nanoparticles and molecules.

It is possible to find field of fluid velocities all over the filter by numerical solving of hydrodynamic equations but it takes a lot of time. It is much better to derive approximate analytical formula for flow through the cell and reduce the problem to the system of algebraic equations. Let us use the formula from [17] which is worthy for cells and apertures of arbitrary form. Flow through the aperture is approximately equals

$$F = -0.8 \frac{p' S^2}{P^2 \mu} s, \qquad (20)$$

where $p'$ – pressure derivative, $S$ – area of the cell section near the aperture, $P$ – perimeter of that cell section, $s$ – area of the aperture, $\mu$ – dynamical viscosity. It is supposed in the formula that liquid flows as if there were no partition, and the partition with the aperture plays the role of a flow limiter. Analog of hydraulic radius is also used: $R = 2S/P$.

Let pressures in the centers of the cells be variables. Pressure derivative $p'$ for the formula (20) is approximately difference of pressures in the centers divided by distance between the centers. So, we obtain the system of linear equations to find pressures:

$$F_{x(i,j,k)} + F_{(i,j,k)x} + F_{y(i,j,k)} + F_{(i,j,k)y} + F_{z(i,j,k)} + F_{(i,j,k)z} = 0,$$
$$(i,j,k) \in U_{вн}, \qquad (21)$$

where $U_{in}$ – set of indexes of the inner cells (it is all the cells excepting ones in inlet and outlet). The sum of flows through all the apertures of the cell is in the left-hand side. The flow is zero if the facet is a part of the filter cover, or some particle has blocked the aperture, or sediment has completely filled the aperture. It is possible to solve the system (21) by Seidel iteration method. The system can become degenerate after closing of too many apertures.

**Some useful formulas.**

It was already said that it is better if liquid flows free in the directions $Ox$ and $Oy$, i.e. if the filter is not porous but a sequence of layers (membranes). Thus, let apertures in directions $Ox, Oy$ be maximally large circles and not filtering. Then flow through the filter is approximately

$$F \approx \frac{1}{2}F_1 \approx \frac{-\pi\, h^2 \langle r \rangle^2 p'}{40\, \mu} n_x n_y, \qquad (22)$$

where $h$ – length of the edge of the cubic cell, $\langle r \rangle$ – average radius of the filtering apertures, $F_1$ – flow throw the central cell of the uncontaminated filter. If it is accepted that the filter can work until more than a half filtering apertures open then

$$FNT \approx \frac{1}{2}F_1 NT n_x n_y \approx \frac{1}{2}(n_x - 1)(n_y - 1)(n_z - 1), \qquad (23)$$

where $T$ – approximate time of the filter work neglecting sedimentation, $N$ – concentration of particles. A useful formula for calculation of the filter follows from (22), (23):

$$NT \approx \frac{-20\, \mu\, n_z}{\pi\, h^2 \langle r \rangle^2 p'} \approx \frac{n_z}{F_1}. \qquad (24)$$

It is supposed here that probability for the particle to flow through all the filter

$$q^{n_z - 1} \qquad (25)$$

is small. The formula (24) gives very large time of the filter work. Real time is much less because the first membrane proves to be totally contaminated earlier. Taking into account sedimentation decreases the time of work additionally. It seems to be a fault of the filter construction and not a fault of the formula. Some variants to solve this problem are at the end of this article.

Concentrations in the layers of the cells are

$$N_k = Nq^{k-1}, k = 1, 2, \ldots, N_z. \qquad (26)$$

where $k$ – number of the layer of cells along axis $Oz$. Taking into account sedimentation it is necessary to calculate average probabilities of penetration through the layers in (26) instead of

$q^{k-1}$.

Calculation process of the filter work at each time step $\Delta t$ is following: calculation of pressures, firstly, then blocking of some apertures by means of random-number generator (filtration), then decreasing of the radiuses of the other apertures (sedimentation). The value $\Delta t$ determines width of sediment (10) and probability of blocking

$$q_{\Delta t} = 1 - q^{F\Delta tN}, \qquad (27)$$

where $F$ – flow through the aperture, $N$ – concentration before the aperture.

Concerning (27) there is some contradiction. Statistically, velocity of contamination is

$$N<(1-q)F>n, \qquad (28)$$

where $<>$ means average over the layer, n – number of open apertures in the layer. I.e. approximately

$$n\left(1 - e^{-N<(1-q)F>\Delta t}\right), \qquad (29)$$

apertures will be blocked during $\Delta t$. And the formula (27) provides that approximately

$$n<1 - q^{F\Delta tN}>, \qquad (30)$$

apertures will be blocked at the same time. So, the probability

$$q_{\Delta t} = (1 - q^{F\Delta tN})\frac{1 - e^{-N<(1-q)F>\Delta t}}{<1 - q^{F\Delta tN}>}, \qquad (31)$$

gives more precious velocity of contamination.

**Results of numerical calculations.**

Parameters of the filter: $L_x = L_y = L_z = 10^{-3}\ m$, $n_x = n_y = n_z = 20$, so, $h_x = h_y = h_z = 5 \cdot 10^{-5}\ m$. Pressure gradient of the filter $p' = -10^4\ Pa\ m^{-1}$. Radiuses of the filtering apertures $r = 1.19 \cdot 10^{-5}\ m$, radiuses of the side apertures $r = 2.5 \cdot 10^{-5}\ m$, length of the particles $l = 2.5 \cdot 10^{-5}\ m$, $q^{n_z-1} = 0.001$ (25), concentration of the particles $1.389 \cdot 10^7\ m^{-3}$. Inlet and outlet are from 6 to 14 cells at both axes. Parameters for sediment are in the item "A variant of chemical reaction".

The filter worked a little more than 3 days until complete stopping of the flow by sediments. Speed near the cover is slow firstly but when the first membranes become blinded in the middle, the flow shifted to the outer boundaries of the filter. Filtration took place with decrease of the radiuses of the apertures. Then began fast disappearing of the apertures because of sediments. Fig. 2 represents the state of the first and the tenth membranes just before the apertures begin to disappear. $v_0 > v_{stat}$ in all open apertures in this case. Up to this moment 271 apertures caught the particles. The filter proved to be low-effective in result: approximately 7300 filtering apertures were closed by sediment. The first membrane catch many particles and the tenth – only one. Line 1 on Fig. 3 shows contamination of the membranes. There are almost no caught particles beginning from 5th membrane but deleting of these membranes will deteriorate pureness of filtration.

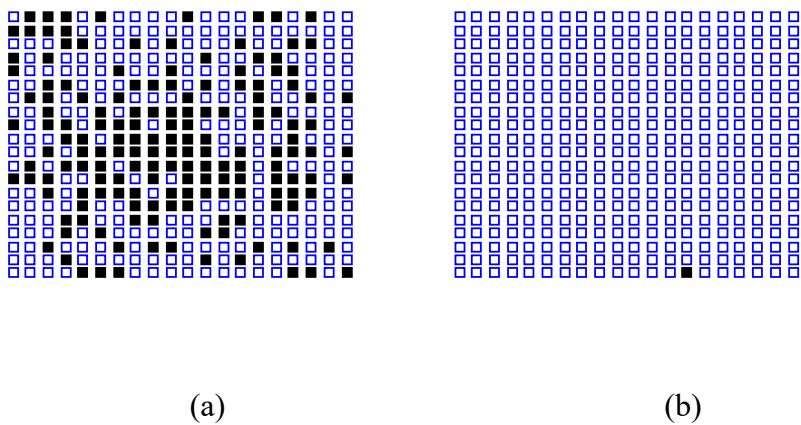

(a)     (b)

*Fig. 2 (a) The first membrane, (b) the tenth membrane. Time 3 days. Empty square – there is flow through the aperture, filled square – the aperture is blocked by the particle.*

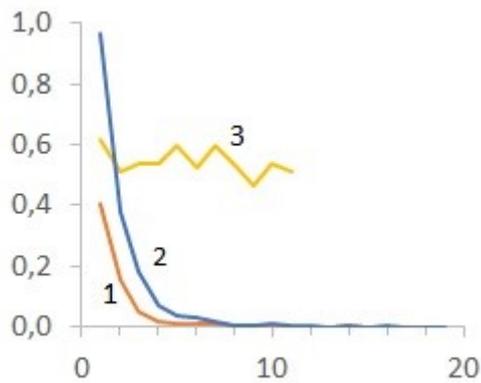

*Fig. 3 Part of contaminated apertures in the membranes. Horizontal axis – number of the membrane. Line 1 – calculation represented on Fig. 2, line 2 – Fig. 4 (increased concentration of the particles), line 3 – Fig. 5 (calculation for equal contamination of all the membranes).*

Two reasons cause ineffective usage of the filtering apertures. The particles rarely achieve the far membranes because of concentration decrease, firstly. Secondly, real time of work is similar to time of blinding of the apertures by sediments.

Fig. 4 represents the first and the tenth membranes in the case when calculated time of work is decreased in 3 times. The first membrane is completely blocked partially by the particles and partially by sediments. The tenth membrane caught much more particles than in previous case. It is also completely blocked, but mainly by sediments. The quantity of the blocked apertures in the filter is 691 at this time and 6440 are closed by sediments. Liquid already do not flow through the filter, so, the other apertures will not catch particles. But comparing with the previous case usage of the filtering apertures increased in two and a half times. Line 2 on Fig. 3 shows contamination of the membranes. There are almost no caught particles beginning from 11th membrane.

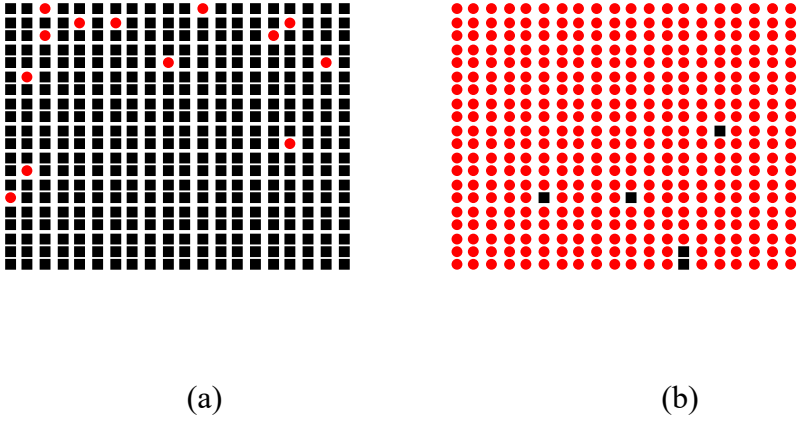

(a)                                          (b)

*Fig. 4 (a) The first membrane, (b) the tenth membrane. Time 3.25 days. Filled circle – the aperture is closed by sediment.* $N = 4.167 \cdot 10^7 \ m^{-3}$.

The one more way to make such a filter more efficient is increasing of the liquid speed. It reduces calculated time of the filter work. However, increasing of speed can cause vortexes. It can also increase velocity of sediment growth (10), (11). There were no notable acceleration of sediment growth in this example, but there are still open apertures in the filter. That is because of dependence of velocity of sedimentation on speed of the flow.

**Different aperture radiuses from layer to layer.**

It is possible to raise a deal of apertures that catch particles by changing apertire radiuses from layer to layer. One way is to make equal the velocity of the layers contamination. If the layers contaminated equally then velocities of their contamination are proportional to $N_k q'_k$, where $N_k$ concentration before the layer, $q'_k = 1 - q_k$ is probability of particle catching and $k$ is a number of the layer. Here $q_k$ is probability for the particle to pass through the $k$-th membrane under the condition that it is already in the $k$-th layer of the cells. To continue the equality of contamination it is necessary to satisfy equalities

$$N_{k+1}q'_{k+1} = N_k q'_k, k = 1,2,\ldots,n_z - 2, \qquad (32)$$

but $N_{k+1} = N_k(1 - q'_k)$, so, the following equation determines the probabilities

$$q'_{k+1} = \frac{q'_k}{1 - q'_k}, k = 1,2, \ldots, n_z - 2. \tag{33}$$

This recurrent formula can be transformed to the usual formula

$$q'_k = \frac{q'_1}{1 - (k-1)q'_1}, k = 2,3, \ldots, n_z - 1. \tag{34}$$

The necessary and sufficient condition for all $q'_k$ be in interval $(0,1)$ is $q'_1 < 1/(n_z - 1)$.

Probability of penetration of the particle throw the filter is

$$\prod_{1}^{n_z-1} 1 - q'_k = 1 - (n_z - 1)q'_1. \tag{35}$$

If $q'_1 = 1/(n_z - 1)$ then $q'_{n_z-1} = 0$, thus, it is always possible to choose $q'_1$ so close to $1/(n_z - 1)$ that probability of particle penetration throw the filter will be enough close to zero. Formula (19) provides a way to obtain necessary $q'_k$ by making apertures of proper radiuses:

$$r_k = \frac{l}{2}\sqrt{1 - (q'_k)^2}. \tag{36}$$

Such a filter has a peculiarity: if the number of membranes is small and pureness of filtration is high, then the apertures in the last membrane an order of magnitude less than in other membranes. That is because function (34) becomes unit at $k = n_z - 1$ and $q'_1 = 1/(n_z - 1)$. All the membranes are really contaminated equally, but liquid flow proves to be very slow. To process this situation it is possible to make many additional apertures of the same radius in the last membrane. As a result, pureness of filtration will be the same and the flow will be faster. Number of particles in all the membranes will be approximately equal, but the ratio of particles to the number of apertures will be less in the last membrane.

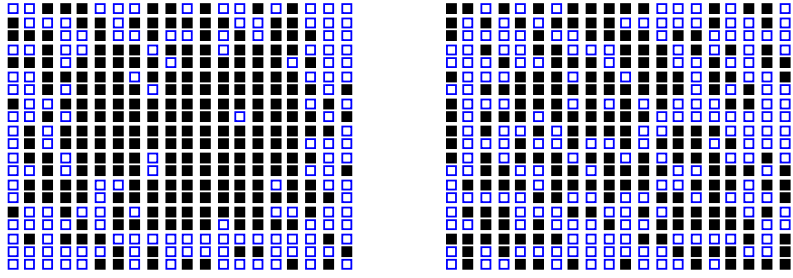

*Fig. 5. The first and the last (11ᵗʰ) membranes. Time 3 days.* $N = 6.165 \cdot 10^7 \ m^{-3}$.

Line 3 on Fig. 3 shows contamination of the membranes in the following calculation: 11 membranes arranged according to (34), probability of penetration through the filter 0.01, calculated time of the filter work 2.5 days, the other parameters from the second numerical calculation neglecting sedimentation. More than a half of apertures have caught particles. The further calculation shows that at least 3968 apertures of 4400 will catch particles. The calculated time of the filter work is really similar to the time when a half of apertures have caught particles.

There is also another interesting way to set radiuses of apertures in the layers. The formula (19) is a distribution function for the value of the half of projection length. Dividing the interval of probabilities (0,1) in $n_z$ equal intervals we can set probabilities for membranes

$$q'_k = \frac{k}{n_z}, k = 1, 2, \dots, n_z - 1. \tag{37}$$

Probability of penetration of the particle throw the filter is

$$\frac{(n_z - 1)!}{(n_z - 1)^{(n_z - 2)}}. \tag{38}$$

To provide pureness of filtration near 0.001 (like in the first numerical example) it is necessary 11 or 12 membranes. I.e. the filter will be almost two times shorter than in the first and the second numerical examples, and its usage of apertures will be much better.

**Conclusion.**

The mathematical model of sediment growth in the filter is constructed in the article and the model of the filter is improved. Analysis of equations and results of calculation lead to the following conclusions.

Usage of the filtering apertures increases if time of the filter work is less than time until blinding of the apertures by sediment. It is possible to reduce the calculated time of the filter work by varying of parameters in the formula (24). Changing of pressure gradient or radiuses of the filtering apertures can increase usage of the apertures in times. It is also necessary to take into account change of velocity of sediment forming. It is rational to dispense flow all over the first membrane. Differing of radiuses of apertures from membrane to membrane can make contamination equal along the filter. Preliminary calculation of the filter for the purpose it will serve can make it much more efficient.

**Literature.**